\documentclass[12pt,a4paper,reqno]{amsart}

\usepackage{anysize}
\usepackage[latin1]{inputenc}
\usepackage{amsmath}
\usepackage{amssymb}
\usepackage{array}

\marginsize{2.5cm}{2.5cm}{3.5cm}{3.5cm}

\begin{document}

\newtheorem{observacion}{Remark}
\newtheorem{theorem}{Theorem}[section]
\newtheorem{corollary}[theorem]{Corollary}
\newtheorem{lemma}[theorem]{Lemma}
\newtheorem{proposition}[theorem]{Proposition}
\newtheorem{definition}[theorem]{Definition}
\newtheorem{example}[theorem]{Example}
\newtheorem{remark}[theorem]{Remark}
\newtheorem{notation}[theorem]{Notation}
\newtheorem{question}[theorem]{Question}
\newtheorem{claim}[theorem]{Claim}
\newtheorem{conjecture}[theorem]{Conjecture}

\title{Kind of proofs of Ramanujan-like series}
\author{Jesús Guillera}

\date{}

\begin{abstract}
We make a summary of the different types of proofs adding some new ideas. In addition we conjecture some relations which could be necessary in ``modular type proofs" (not still found) of the Ramanujan-like series for $1/\pi^2$.
\end{abstract}

\maketitle

\section{Ramanujan-type series for $1/\pi$}

\subsection{Introduction}

We recall that a Ramanujan-type series for $1/\pi$ is a series of the form
\[ \sum_{n=0}^{\infty} \frac{ \left( \frac{1}{2} \right)_n (s)_n (1-s)_n}{ (1)_n^3} z^n (a(z)+b(z)n)= \frac{1}{\pi}, \]
where $s=1/2$, $1/4$, $1/3$, or $1/6$ and $z$, $a(z)$, $b(z)$ are algebraic numbers. The name is in honor to Ramanujan who discovered $17$ series of this form in 1914. One of them is the celebrated formula
\begin{equation}\label{rama42}
\sum_{n=0}^{\infty} \frac{1}{2^{6n}} \; \frac{\left( \frac12 \right)_n^3}{(1)_n^3} \; (42n+5)=\frac{16}{\pi}.
\end{equation}
The brothers Jonathan and Peter Borwein were the first (1987) to prove the $17$ series discovered by S. Ramanujan. An excellent survey is \cite{BaBeCh}.

\subsection{q-parametrization}
The function
\[ w_0=\sum_{n=0}^{\infty} \frac{ \left( \frac{1}{2} \right)_n (s)_n (1-s)_n}{ (1)_n^3} z^n, \]
satisfies the differential equation
\[ \left( \theta^3-z(\theta+1/2)(\theta+s)(\theta+1-s) \right)w=0, \quad \theta=z \frac{d}{dz}. \]
This equation has three fundamental solutions $w_0$, $w_1$, $w_2$. As usual we let
\[ q=\exp (\frac{w_1}{w_0}). \]
Writing $q$ as a series of powers of $z$ and inverting it we have $z$ as a series of powers of $q$. As $z$ is a function of $q$ then we can also view $b$ and $a$ as functions of $q$. It is known that $z(q)$ is a modular function and that $b(q)$ and $a(q)$ are modular forms of weight $2$. Let $q=e^{i \pi r} e^{-\pi \tau}$, with $\tau>0$. The interesting fact is that we have Ramanujan-type series if and only if $r$ and $\tau^2$ are rational. We will only consider the cases $r=0$ (series of positive terms), that is $q=e^{-\pi \tau}$ and $r=1$ (alternating series), that is $q=-e^{-\pi \tau}$. The other values of $r$ lead to complex series \cite{AlGu2}. For a fix value of $r$ we can consider $z$, $b$, $a$ as functions of $\tau$.

\subsection{Modular equations}
A modular equation of order $k$ for a modular function $z(q)$ is an expression relating $z(q)$ with $z(q^k)$ in an algebraic way. In the examples we show how to determine $z$ combining the functional relation (it depends only on $s$) with a modular equation. Then we determine $b$ using a simple known formula. The value of $a$ is more difficult to determine.

\subsection*{Example 1}

Let $q=e^{-\pi \tau}$ with $\tau>0$. For $s=1/2$ it is known that
\[ z(\tau)=4 \lambda(\tau)(1-\lambda(\tau)), \]
where $\lambda(\tau)$ is a modular function which is known as the elliptic $\lambda$-function. We see that (\ref{rama42}) has $z=1/2^6$. To understand the origin of this value we substitute $\tau=\sqrt{7}$ in the functional relation $1-\lambda(\tau)=\lambda(1/\tau)$ and in the septic modular equation
\[
\left\{ \frac{}{} \! \! \lambda(\tau)\lambda(\tau/7) \right\}^{\frac18}+\left\{ \frac{}{} \! \! (1-\lambda(\tau))(1-\lambda(\tau/7)) \right\}^{\frac18}=1.
\]
This septic modular equation is due to C. Guetzlaff $(1834)$ and was rediscovered by Ramanujan. Then to find the value of $b$ we can use the simple formula $b(\tau)=\tau \sqrt{1-z(\tau)}$ \cite{Gu5}. The value of $a$ is more difficult to prove
\cite[Th. 9.1]{BaBeCh}.

\subsection*{Example 2}

For the alternating series corresponding to $s=1/2$ we can use the known duality identity $z(\tau) \, z(2/\tau)=1$ as we do in this example. Let $q=-e^{-\pi \tau}$ with $\tau>0$. We observe that Ramanujan-type series
\[ \sum_{n=0}^{\infty} \frac{(-1)^n}{2^{3n}} \frac{\left( \frac12 \right)_n^3}{(1)_n^3} (6n+1)=\frac{2 \sqrt{2}}{\pi}, \]
has $z=-1/8$. To understand the origin of this value we use the duality identity and the modular equation
\[ 64x+y^2-48xy+64x^2y=0, \]
where $x=z(\tau)$ and $y=z(\tau/2)$. For $\tau=2$, we get $128x^3-48x^2+1=0$. One solution is $x=-1/8$. Then, from the formula
$b=\tau \sqrt{1-z}$, we obtain $b=3/\sqrt2$.

\subsection*{Example 3}

It is known that there are modular equations of any order. For example, to find a modular equation of third order for the case $s=1/2$, it is enough to consider
\[ P(x,y)=\sum_{i=0}^6 \sum_{j=0}^{6-i} a_{ij} x^i y^j = 0, \]
where $x=z(q)$, $y=z(q^3)$. Taking $50$ terms in $z(q)$ and in $z(q^3)$, and solving the linear system of equations formed by the coefficients of the powers of $q$, we guess that
\begin{multline}\nonumber
-4096xy + 4608(x^2y+xy^2)+(x^4+y^4)-900(x^3y+xy^3) \\ +28422x^2y^2 + 4608(x^3y^2+x^2y^3) - 4096 x^3y^3 = 0.
\end{multline}
A much more difficult task is to prove the modular equations rigourously. H.H. Chan hints in \cite{Chan2} that the above modular equation can be proved using \cite[p. 231, xii]{Be3}.

\subsection{Hypergeometric transformations}

The following hypergeometric identities:
\begin{align}
&\sum_{n=0}^{\infty} \frac{\left( \frac12 \right)_n^3}{(1)_n^3} z^n=\frac{1}{\sqrt{1-z}} \sum_{n=0}^{\infty} \frac{\left( \frac12 \right)_n\left( \frac14 \right)_n\left( \frac34 \right)_n}{(1)_n^3} \left( \frac{-4z}{(1-z)^2} \right)^n
\nonumber \\ &=\frac{2}{\sqrt{4-z}} \sum_{n=0}^{\infty} \frac{\left( \frac12 \right)_n\left( \frac16 \right)_n\left( \frac56 \right)_n}{(1)_n^3} \left( \frac{27z^2}{(4-z)^3} \right)^n
\! = \! \frac{1}{\sqrt{1-4z}} \sum_{n=0}^{\infty} \frac{\left( \frac12 \right)_n\left( \frac16 \right)_n\left( \frac56 \right)_n}{(1)_n^3} \left( \frac{-27z}{(1-4z)^3} \right)^n, \nonumber
\end{align}
see \cite[eqs. 3.3, 3.7--3.9]{BaBe},
\begin{equation}\label{rogers}
\sum_{n=0}^{\infty} \frac{\left( \frac12 \right)_n \left( \frac14 \right)_n \left( \frac34 \right)_n }{(1)_n^3} \left( \frac{256z^3}{9(3+z)^4} \right)^n = \frac{3+z}{3(1+3z)} \sum_{n=0}^{\infty} \frac{\left( \frac12 \right)_n \left( \frac14 \right)_n \left( \frac34 \right)_n}{(1)_n} \left(\frac{256z}{9(1+3z)^4}\right)^n,
\end{equation}
\cite[eq. 3.7]{Ro}, and also others of the same style like \cite[eqs. 3.3--3.6]{BaBe} and \cite[eq. 3.6]{Ro} are known. They can be proved by purely hypergeometric methods.

\subsection{Zudilin's translation method}

We explain it with two examples. The method is in fact much more powerful; for example, it is applied in \cite{WaZu} in much more general settings. We can translate Ramanujan-type real series into Ramanujan-type real series (that is series with $q=e^{-\pi \tau}, \tau>0$ or with $q=-e^{-\pi \tau}, \tau>0$) if $\tau_2/\tau_1$ is rational (otherwise it is not possible). W. Zudilin shows how to find the hypergeometric transformation one needs to do it. Once we have found it, he shows that we can prove it without modularity considerations \cite{Zu5}.

\subsection*{Example 1}

We have
\[
\frac{1}{\pi}=\frac{\sqrt2}{4} \sum_{n=0}^{\infty} \frac{\left( \frac12 \right)_n^3}{(1)_n^3} \left( \frac{-1}{8} \right)^n (6n+1)=\frac{2}{11 \sqrt{33}} \sum_{n=0}^{\infty} \frac{\left( \frac12 \right)_n
\left( \frac16 \right)_n \left( \frac56 \right)_n}{(1)_n^3} \left( \frac{2}{11} \right)^{3n} (126n+10).
\]
The first series was proved by the WZ-method (see \cite{Gu-on-wz} and its references). We prove the other one translating it in the following way
\[
\left( 1 + 6 z \frac{d}{dz} \right) \left\{ \sum_{n=0}^{\infty} \frac{\left( \frac12 \right)_n^3}{(1)_n^3} z^n \right\} = \left( 1 + 6 z \frac{d}{dz} \right) \left\{ \frac{2}{\sqrt{4-z}} \sum_{n=0}^{\infty} \frac{\left( \frac12 \right)_n\left( \frac16 \right)_n\left( \frac56 \right)_n}{(1)_n^3} \left( \frac{27z^2}{(4-z)^3} \right)^n \right\}.
\]
Finally we take $z=-1/8$. The original proof of the series in the right side is due to the Borweins \cite{Bo2}, and is based in the modular theory.

\subsection*{Example 2}

We have
\begin{align}
\frac{1}{\pi} &=\frac{\sqrt{3}}{16} \sum_{n=0}^{\infty} \frac{\left( \frac12 \right)_n \left( \frac14 \right)_n \left( \frac34 \right)_n}{(1)_n^3} \left( \frac{-1}{48} \right)^n (28n+3) \nonumber \\
&= \frac{3}{\sqrt[4]{12}} \sum_{n=0}^{\infty} \frac{\left( \frac12 \right)_n^3}{(1)_n^3} (2-\sqrt{3})^{4n} \left( (8\sqrt3 - 12)n+(3\sqrt3-5) \right). \nonumber
\end{align}
The first series was proved by the WZ-method (see \cite{Gu-on-wz} and its references). We prove the other one translating it in the following way
\[
\left( 3+28 u \frac{d}{du} \right) \left\{ \sum_{n=0}^{\infty} \frac{\left( \frac12 \right)_n \left( \frac14 \right)_n \left( \frac34 \right)_n}{(1)_n^3} u^n \right\} = \left( 3 + \frac{28 u}{z u'} z \frac{d}{dz} \right) \left\{ \sqrt{1-z} \sum_{n=0}^{\infty} \frac{\left( \frac12 \right)_n^3}{(1)_n^3} z^n \right\},
\]
where $u=(-4z)/(1-z)^2$ and $u'=du/dz$. Finally we substitute $u=-1/48$. The first series is due to Ramanujan, the other one is due to the Borweins \cite{Bo}, and the original proofs are based in the modular theory.

\subsection{A variant of Zudilin's method}

Alexander Aycock (a student of the Johannes-Gutenberg-Universität Mainz), had the idea of applying the translation to a kind of limit cases. This example is essentially due to him. Let $u=256z^3/(9(3+z)^4)$ and $u'=du/dz$. Applying to (\ref{rogers}) the operator $u \cdot d/du$, we have
\[
\sum_{n=0}^{\infty} \frac{\left( \frac12 \right)_n \left( \frac14 \right)_n \left( \frac34 \right)_n }{(1)_n^3} n u^n = \frac{u}{u'} \frac{d}{dz} \left[ \frac{3+z}{3(1+3z)} \sum_{n=0}^{\infty} \frac{\left( \frac12 \right)_n \left( \frac14 \right)_n \left( \frac34 \right)_n}{(1)_n} \left(\frac{256z}{9(1+3z)^4}\right)^n \right].
\]
Then we multiply by $\sqrt{1-u}$ and take the limit as $u \to 1^{-}$. Writing $(1-u)^{-1/2}$ as a series of powers of $u$ and applying  the Stolz-Cesàro theorem we see that the limit of the left side is equal to $\sqrt{2}/(2 \pi)$. Hence we have (observe that $z$ tends to $9^{+}$)
\begin{equation}\label{aycock}
\lim_{z \to 9^{+}} \left\{ \sqrt{1-u} \frac{u}{u'} \frac{d}{dz} \left[ \frac{3+z}{3(1+3z)} \sum_{n=0}^{\infty} \frac{\left( \frac12 \right)_n \left( \frac14 \right)_n \left( \frac34 \right)_n}{(1)_n} \left(\frac{256z}{9(1+3z)^4}\right)^n \right] \right\}=
\frac{\sqrt{2}}{2 \pi},
\end{equation}
Making the calculations it yields
\begin{equation}\nonumber
\sum_{n=0}^{\infty} \frac{1}{7^{4n}} \; \frac{\left( \frac12 \right)_n  \left( \frac14 \right)_n \left( \frac34 \right)_n}{(1)_n^3} \; (40n+3) = \frac{\sqrt3}{9\pi},
\end{equation}
which is one of the series for $1/\pi$ discovered by Ramanujan \cite[eq. 42]{Ra}.

\subsection{Non-modular-type proofs versus modular-type proofs}

Other kind of hypergeometric proofs for some Ramanujan-type series are given by W. Chu in \cite{Chu} by acceleration of Dougall's bilateral $_2H_2$ series. If we use Zudilin's method to translate series proved by the WZ-method (see \cite{Gu-on-wz} and its references) then the proofs are purely hypergeometric. Now, see the formulas of the papers \cite{BaBe} and \cite{BaBe2}, where $n$ in \cite{BaBe} and $2n$ in \cite{BaBe2} stand respectively for $\tau^2$. As we have proved by the WZ-method a series (at least) in each of the cases $\tau^2=2,3,4,5,6,7,9$, we can prove by ``translation" all the formulas corresponding to $\tau^2=2,3,4,5,6,7,9$ and also those for $\tau^2=18$ and $\tau^2=25$. In addition, we can derive as well the formulas with $\tau^2=15$ translating the following ``divergent" series for $1/\pi$:
\[
\sum_{n=0}^{\infty} \frac{\left( \frac12 \right)_n  \left( \frac13 \right)_n \left( \frac23 \right)_n}{(1)_n^3} \; (-4)^n (15n+4) \, \text{``=''} \, \frac{3\sqrt3}{\pi},
\]
which has $\tau=b/\sqrt{1-z}=\sqrt{15}/3$ and was proved in \cite{Gu8} by the WZ-method. However we cannot prove those formulas with $\tau^2=10,11,13,14,17,22$. The unique known proofs for them and others even more complicated is by using the modular theory, as Nayandeep Baruah and Bruce Berndt do in \cite{BaBe} and \cite{BaBe2}. In $1987$ J. and P. Borwein gave modular-type proofs of the two most impressive series discovered by Ramanujan ($\tau^2=35$ and $\tau^2=58$ respectively), namely:
\begin{align}
\sum_{n=0}^{\infty} \frac{(-1)^n}{882^{2n}} \; \frac{\left( \frac12 \right)_n  \left( \frac14 \right)_n \left( \frac34 \right)_n}{(1)_n^3} \; (21460n+1123) &= \frac{3528}{\pi}, \nonumber \\
\sum_{n=0}^{\infty} \frac{1}{99^{4n}} \; \frac{\left( \frac12 \right)_n  \left( \frac14 \right)_n \left( \frac34 \right)_n}
{(1)_n^3} \; (26390n+1103) &= \frac{9801 \sqrt 2}{4 \pi}, \nonumber
\end{align}
which give almost $6$ and $8$ digits per term respectively \cite{Bo}. The scope of of the WZ-method is unknown. For example, we ignore if these two last series or simpler series, like ($\tau^2=10$):
\begin{equation}\label{rama-10-1}
\sum_{n=0}^{\infty} \frac{1}{3^{4n}} \; \frac{\left( \frac12 \right)_n  \left( \frac14 \right)_n \left( \frac34 \right)_n}{(1)_n^3} \; (10n+1) = \frac{9\sqrt2}{4 \pi},
\end{equation}
\cite[eq. 41]{Ra} are provable by this method.

\subsection*{Remark}
We propose a way of proving (\ref{rama-10-1}) without using modular equations, based in a generalization of Aycock's idea. Here we only sketch the proof. Begin with the following limit:
\begin{equation}\label{limit-yy}
\lim_{z \to 1^{-}} \sqrt{(1-z)\left(1+\frac{z}{4}\right)} \sum_{n=0}^{\infty} \sum_{k=0}^{n} \binom{n}{k}^4 n \left( \frac{z}{16} \right)^n = \frac{1}{\pi \tau_c} = \frac{\sqrt{10}}{2 \pi},
\end{equation}
where $q_c=e^{-\pi \tau_c}$ and $z(q_c)=1$. As $\tau_c=2/\sqrt{10}$ and the value of $\tau$ corresponding to (\ref{rama-10-1}) is $\tau=5 \tau_c$, we can find a transformation \cite{Zu5} which allows to translate (\ref{limit-yy}) into the series (\ref{rama-10-1}). Our example is a particular case of the general limit
\begin{equation}\label{limit-general}
\lim_{z \to z_c^{-}} \sqrt{P(z)} \sum_{n=0}^{\infty} A_n n z^n = \frac{1}{\pi \tau_c},
\end{equation}
where $A_n$ are a type of Ramanujan-Sato numbers, $P(z)$ is the polynomial defined in \cite{AlGu2}, $z_c$ is the radius of convergence of the series, and $\tau_c$ is defined implicitly by $q_c=e^{-\pi \tau_c}$ and $z_c=z(q_c)$. The proof of (\ref{limit-general}) uses the formula $b=\tau \sqrt{P(z)}$ obtained in \cite{AlGu2}. We thank to Anton Mellit the idea of applying Stolz-Cesáro theorem to get the limit in (\ref{aycock}), and in general the limits of the form (\ref{limit-general}) assuming we know the asymptotic behavior of $A_n$ \cite{McI}. This method allows to determine the exact value of $\tau_c$.

\section{Ramanujan-like series for $1/\pi^2$ (A summary)}

\subsection{Introduction}
Let $s_0=1/2 \,$, $s_3=1-s_1 \,$, $s_4=1-s_2 \,$. We recall that a Ramanujan-like series for $1/\pi^2$ is a series of the form
\[
\sum_{n=0}^{\infty} z^n \left[\prod_{i=0}^{4} \frac{(s_i)_n}{(1)_n}\right] (a+bn+cn^2)=\frac{1}{\pi^2},
\]
where $z$, $a$, $b$ and $c$ are algebraic numbers and the possible couples $(s_1, s_2)$ are $(1/2,1/2)$, $(1/2,1/3)$, $(1/2,1/4)$, $(1/2,1/6)$, $(1/3,1/3)$, $(1/3,1/4)$, $(1/3,1/6)$, $(1/4,1/4)$, $(1/4,1/6)$, $(1/6,1/6)$, $(1/5,2/5)$, $(1/8,3/8)$, $(1/10,3/10)$, $(1/12,5/12)$. Up till now only the following series have been proved:
\begin{align}
\sum_{n=0}^{\infty} \frac{\left( \frac12 \right)_n^3 \left( \frac14 \right)_n \left( \frac34 \right)_n}{(1)_n^5} \frac{1}{2^{4n}} (120n^2+34n+3) &=\frac{32}{\pi^2}, \label{wz1-pi2} \\
\sum_{n=0}^{\infty}\frac{\left(\frac12 \right)_n^5}{(1)_n^5}\frac{(-1)^n}{2^{10n}} (820n^2+180n+13) &=\frac{128}{\pi^2}, \label{wz2-pi2} \\
\sum_{n=0}^{\infty} \frac{\left( \frac12 \right)_n^5}{(1)_n^5} \frac{(-1)^n}{2^{2n}} (20n^2+8n+1) &=\frac{8}{\pi^2}, \label{wz3-pi2} \\
\sum_{n=0}^{\infty}  \frac{\left(\frac12 \right)_n^3 \left( \frac13 \right)_n \left( \frac23 \right)_n}{(1)_n^5} \left( \frac34 \right)^{3n} (74n^2+27n+3) &=\frac{48}{\pi^2}. \label{wz4-pi2}
\end{align}
All the known proofs are of hypergeometric type. In $2002$ we proved (\ref{wz1-pi2}) and (\ref{wz2-pi2}) by the WZ-method and in $2003$ and $2010$, again by the WZ-method, we proved respectively (\ref{wz3-pi2}) and (\ref{wz4-pi2}), (see \cite{Gu-new} and its references). In \cite{Gu7} we show a collection of conjectured series for $1/\pi^2$. In $2010$, in a joint paper with Gert Almkvist \cite{AlGu}, a new series with $(s_1,s_2)=(1/3,1/6)$ was conjectured. W. Zudilin was the first to realize that this kind of series were related to the theory of Calabi-Yau threefolds \cite{YaZu} and \cite{Zu}. We believe that this is the start point towards proofs based in ``modularity" \cite{Zu4} and \cite{Mo}.

\subsection{Calabi-Yau differential equation}

The hypergeometric function
\[ w_0=\sum_{n=0}^{\infty } \left[ \prod_{i=0}^{4} \frac{(s_i)_n}{(1)_n}\right] z^n, \]
satisfies a fifth order differential equation
\[ \left( \theta^5-z \prod_{i=0}^4 (\theta+s_i) \right) w_0=0. \]
This differential equation is of a special type. Its $5$ fundamental solutions can be recovered from the $4$ solutions $y_{0},y_{1},y_{2},y_{3}$ of a Calabi-Yau differential equation. The parametrization
\[ q=\exp(\frac{y_1}{y_0}) \]
defines $z$ as a series of powers of $q$. This function $z(q)$ is called the mirror map. The Yukawa coupling is then defined as
\[ K(q)=\theta_q^2 (\frac{y_2}{y_0}), \qquad \theta_q=q \frac{d}{dq}. \]
The unique power series $T(q)$, such that $T(0)=0$, which is related to the Yukawa coupling in the way
\[ \theta^3_q \, T(q)=1-K(q), \]
plays an important role in the theory. From it we define the crucial functions:
\begin{equation}\label{alfa-tau}
\alpha(q)=\frac{\frac{1}{6} \ln^3|q|-T(q)-h \zeta (3)}{\pi^2 \ln|q|}, \quad
\tau(q)=\frac{\frac{1}{2} \ln^2|q|-(\theta_q T)(q)}{\pi^2}-\alpha(q),
\end{equation}
where
\begin{equation}\label{h}
h=\frac{2}{\zeta(3)} \left\{ \frac{}{} \zeta(3,1/2)+\zeta(3,s_1)+\zeta(3,s_2)+\zeta(3,1-s_1)+\zeta(3,1-s_2) \right\},
\end{equation}
We consider the values
\begin{equation}\label{criticos}
2 \alpha_c=\frac53+\cot^{2}(\pi s_1)+\cot^{2}(\pi s_2),  \quad \tau_c^2=\frac{1}{\sin^2(\pi s_1) \sin^2(\pi s_2)}.
\end{equation}
They correspond to $q_c$, where $z(q_c)=1$. Sometimes, instead of $\alpha(q)$ we will use the related function
\[ k(q)=2(\alpha(q)-\alpha_c), \]
which translates the critical point to $0$. In \cite{Gu5} we used the functions $k$, $\tau$ and $j$ (which is related to $\tau$ in a simple way); in \cite{AlGu}, we used $k$ and $j$, and in \cite{AlGu2} we considered the functions in (\ref{alfa-tau}) and in addition gave explicit formulas for $c(q)$, $b(q)$ and $a(q)$. The constants (\ref{h}) and (\ref{criticos}) were given in \cite{Gu5} and more explicitly in \cite{AlGu}. For non-hypergeometric functions satisfying Calabi-Yau differential equations see \cite{AlGu} and \cite{AlGu2}. In these cases we cannot (of course) obtain $h$ from (\ref{h}) and $\alpha_c$, $\tau_c$ from (\ref{criticos}) but we can use the method explained in \cite{AlGu2} which consist in finding the critical value $q_c$ solving the equation
\[ \frac{dz}{dq}(q_c)=0, \]
and looking for integer relations among the numbers
\[ \frac{1}{6} \ln^3|q_c| \! - \! T(q_c), \quad \pi^2 \ln|q_c|, \quad \zeta(3). \]
Even better, we can use the explicit formulas given in \cite[Addendum]{AlGu2}.

\subsection{Main conjecture}
The following conjecture was first stated in \cite{Gu5} for hypergeometric series and extended in \cite{AlGu2} to all cases. We use the notation $\alpha_0=\alpha(q_0)$, $\alpha_1=\alpha(q_1)$, $\tau_0=\tau(q_0)$, etc.
\begin{conjecture}\label{conj1}
Let $q=|q|e^{i \pi r}$,
\[ \sum_{n=0}^{\infty} z(q)^n \left[ \prod_{i=0}^4 \frac{(s_i)_n}{(1)_n} \right] (a(q)+b(q)n+c(q)n^2)=\frac{1}{\pi^2}, \]
Then $r$, $\alpha_0$, $\tau_0^2$ are rational if and only if $z_0$, $c_0$, $b_0$, $a_0$ are algebraic.
\end{conjecture}

\section{Ramanujan-like series for $1/\pi^2$ (New content)}

The following conjectures are new:

\begin{conjecture}\label{monodromy}
The monodromy matrix around $z=z_c$ is given by
\[
\frac{1}{\tau_c^2} \left(
  \begin{array}{ccccc}
    \alpha_c^2 & 0 & -\alpha_c(\tau_c^2-\alpha_c^2)/8 & (\tau_c^2-\alpha_c^2)d & -(\tau_c^2-\alpha_c^2)^2/128 \\
    -32d & \tau_c^2 & -8\alpha_c d & 32 d^2 & -(\tau_c^2-\alpha_c^2)d \\
    -8\alpha_c & 0 & \tau_c^2-2\alpha_c^2 & 8 \alpha_c d & -\alpha_c(\tau_c^2-\alpha_c^2)/8 \\
    0 & 0 & 0 & \tau_c^2 & 0 \\
    -32 & 0 & -8 \alpha_c & 32 d & \alpha_c^2 \\
  \end{array}
\right).
\]
\end{conjecture}
For the six examples given in the table of \cite[Th. 3]{ChYaYu} we get the correct matrix. This supports the conjecture.

\begin{conjecture}\label{conj2}
Let $q=e^{-\pi t}$, with $t>0$. If one of the relations $z(t_2)=z(t_1)$ and
\[
\tau_2=\tau_c^2 \cdot \frac{4 \tau_1}{4\tau_1^2-k_1^2}, \qquad k_2=-\tau_c^2 \cdot \frac{4k_1}{4\tau_1^2-k_1^2}, \qquad
t_2=t_1 \cdot \frac{2\tau_1-k_1}{2\tau_1+k_1},
\]
holds all the others are also true.
\end{conjecture}
It is supported by numerical calculations. For example, solving $z(q)=0.9$, we could obtain the identities with a precision of about $25$ digits. To discover them we used the PSLQ algorithm, developed by Ferguson and Bailey, to look for integer relations among the numbers (1): $\tau_1 k_2$ and $\tau_2 k_1$, (2): $k_1 k_2$, $\tau_1 \tau_2$ and $1$, (3): $t_1 t_2^{-1}$, $t_2 t_1^{-1}$, $k_1 k_2$ and $1$. Far from the critical value ($z\!=\!1$ in the hypergeometric cases), for example for $z=0.4$ (hypergeometric cases), we may have divergences but we think the identities still hold by analytic continuation.

\begin{conjecture}(Duality)\label{conj3}
Let $q=-e^{-\pi t}$, with $t>0$ and $(s_1,s_2)=(1/2,1/2)$. If one of the relations $z(t_1)z(t_2)=1$ and
\[
\tau_2=\frac{8 \tau_1}{4\tau_1^2-(k_1+1)^2}, \qquad k_2+1=\frac{8(k_1+1)}{4\tau_1^2-(k_1+1)^2}, \qquad
t_2= \frac{2}{t_1} \cdot \frac{2\tau_1+k_1+1}{2\tau_1-k_1-1},
\]
holds all the others are also true. In that case, we also have the relations
\[
c(t_2)=\frac{\tau_2}{\tau_1} \cdot \frac{c}{\sqrt{-z}}(t_1), \quad
b(t_2)=\frac{\tau_2}{\tau_1} \cdot \frac{c-b}{\sqrt{-z}}(t_1), \quad
a(t_2)=\frac{\tau_2}{\tau_1} \cdot \frac{c-2b+4a}{4 \sqrt{-z}}(t_1).
\]
\end{conjecture}
It is supported by experimental calculations: For example, using the PSLQ algorithm we have found, with a good precision, the following identities:
\begin{align}
(k_1+1)\tau_2 &= (k_2+1)\tau_1, \nonumber \\
(k_1+1)(k_2+1)+8 &= 4\tau_1 \tau_2, \nonumber \\
2 t_1t_2 &= (k_1+1)(k_2+1)+2(k_1+1) \tau_2 + 4. \nonumber
\end{align}
The dual of (\ref{wz3-pi2}) and (\ref{wz2-pi2}) are respectively
\begin{align}
\sum_{n=0}^{\infty} \frac{\left( \frac12 \right)_n^5}{(1)_n^5} (10n^2+6n+1) (-1)^n 2^{2n} & \, \text{``=''} \, \frac{4}{\pi^2}, \nonumber \\
\sum_{n=0}^{\infty} \frac{\left( \frac12 \right)_n^5}{(1)_n^5} (205n^2+160n+32) (-1)^n 2^{10n} & \, \text{``=''} \, \frac{16}{\pi^2}. \nonumber
\end{align}
These two hypergeometric series are ``divergent''. In the following section we will see how to check numerically this kind of series. It is also important to note that W. Zudilin observed a relation among Ramanujan-like series and a kind of p-adic supercongruences \cite{Zu0} and that in \cite{GuiZu} we showed that the same pattern holds when the series is ``divergent". A weaker version of Conj. \ref{conj3} (without the relations for $t$) is in \cite{GuiZu}.

\subsection{Minimal polynomial}

When we discovered and proved the series (\ref{wz4-pi2}) for $1/\pi^2$, we checked that it has $k=2/3$ and $\tau=\sqrt{37}/3$. This was the first known series for $1/\pi^2$ with a non-integer value of $k$. Inspired by it, we tried other values with thirds. For $k=8/3$, using the method explained in \cite{Gu5}, we discovered (also in 2010) the unproved series
\[
\sum_{n=0}^{\infty}  \frac{\left(\frac12 \right)_n^3 \left( \frac13 \right)_n \left( \frac23 \right)_n}{(1)_n^5} (3 \phi)^{3n}
\left[ \! \! \frac{}{}(32-216 \phi)n^2+(18-162\phi)n +(3-30\phi) \frac{}{} \! \! \right] \overset?=\frac{3}{\pi^2},
\]
where
\[ \phi=\left( \frac{\sqrt{5}-1}{2} \right)^5 \simeq 0.09016994374947424102293417182819058860154589902881 \dots. \]
This is the only known convergent series for $1/\pi^2$ with a value of $z$ which is not rational. This value of $z$ has
the minimal polynomial
\[ P(z)=z^2+36828z-729. \]
We suspected that the other root of $P(z)$, which is
\[ z=\left(\frac{-3}{\phi}\right)^3, \]
should be the $z$ of another Ramanujan-like series. To see that we were right we use the analytic continuation given by the hypergeometric series in the right side:
\begin{align}
\sum_{n=0}^{\infty} \frac{\left(\frac12 \right)_n^3 \left( \frac13 \right)_n \left( \frac23 \right)_n}{(1)_n^5} z^n &=
{}_5F_4\left(\begin{matrix}
\frac12, & \frac12, & \frac12, & \frac13, & \frac23, & \nonumber \\[1.1ex]
&  1, & 1, & 1, & 1 \end{matrix} \biggm| z \right), \nonumber \\
\sum_{n=0}^{\infty} \frac{\left(\frac12 \right)_n^3 \left( \frac13 \right)_n \left( \frac23 \right)_n}{(1)_n^5} n z^n &=
\frac{z}{36} \, {}_5F_4\left(\begin{matrix}
\frac32, & \frac32, & \frac32, & \frac43, & \frac53, & \nonumber \\[1.1ex]
&  2, & 2, & 2, & 2 \end{matrix} \biggm| z \right), \nonumber
\end{align}
and
\[
\sum_{n=0}^{\infty} \frac{\left(\frac12 \right)_n^3 \left( \frac13 \right)_n \left( \frac23 \right)_n}{(1)_n^5} n^2 z^n =
\frac{z}{36} \, {}_5F_4\left(\begin{matrix}
\frac32, & \frac32, & \frac32, & \frac43, & \frac53, & \nonumber \\[1.1ex]
&  1, & 2, & 2, & 2 \end{matrix} \biggm| z \right).
\]
Looking for integer relations (with the PSLQ algorithm) among the numbers $v_0$, $v_0 \phi$, $v_1$, $v_1 \phi$, $v_2$, $v_2\phi$ and $1/\pi^2$, where
\[
v_0={}_5F_4\left(\begin{matrix} \frac12, & \frac12, & \frac12, & \frac13, & \frac23, & \nonumber \\[1.1ex]
& 1, & 1, & 1, & 1 \end{matrix} \biggm| \frac{-27}{\phi^3} \right), \quad
v_1=\frac{z}{36} \, {}_5F_4\left(\begin{matrix}
\frac32, & \frac32, & \frac32, & \frac43, & \frac53, & \nonumber \\[1.1ex]
&  2, & 2, & 2, & 2 \end{matrix} \biggm| \frac{-27}{\phi^3} \right),
\]
and
\[
v_2=\frac{z}{36} \, {}_5F_4\left(\begin{matrix} \frac32, & \frac32, & \frac32, & \frac43, & \frac53, & \nonumber \\[1.1ex]
&  1, & 2, & 2, & 2 \end{matrix} \biggm| \frac{-27}{\phi^3} \right),
\]
we find and conjecture the ``divergent" series
\[
\sum_{n=0}^{\infty}  \frac{\left(\frac12 \right)_n^3 \left( \frac13 \right)_n \left( \frac23 \right)_n}{(1)_n^5} \left( \frac{-3}{ \phi} \right)^{3n} \left[ \! \! \frac{}{}(2408+216 \phi)n^2+(1800+162\phi)n+(333+30\phi) \frac{}{} \! \! \right] \text{``\!}\overset?=\text{\!''} \, \frac{36}{\pi^2}.
\]
We can check that the mosaic supercongruences pattern holds for this series \cite{Gu6}. Inspired by these congruences, by the conjectured formulas \cite[Conj. 1.1--1.6]{Sun} and by \cite{GuiRo}, we have guessed that the sum of the series
\[
\sum_{n=0}^{\infty}  \frac{(1)_n^5}{\left(\frac12 \right)_n^3 \left( \frac13 \right)_n \left( \frac23 \right)_n} \left( \frac{-\phi}{ 3} \right)^{3n} \frac{(2408+216 \phi)n^2-(1800+162\phi)n+(333+30\phi)}{n^5}
\]
is equal to
\[ \frac{1125}{4}\sqrt{5}L_5(3)-448\zeta(3). \]
As the convergence of the series is fast, we can use it to get many digits of $L_5(3)$. Another example of the same style is the upside-down series associated to Jim Cullen's formula for $1/\pi^4$ \cite[Sect. 2.5]{Zu4}, namely
\[
\sum_{n=1}^{\infty} \frac{(1)_n^9 2^{12(n-1)}}{\left( \frac12 \right)_n^7 \left( \frac14 \right)_n \left( \frac34 \right)_n} \frac{43680n^4-20632n^3+4340n^2-466n+21}{n^9} \text{``\!}\overset?=\text{\!''} \, -95232 \zeta(5)-160 \pi^5 i.
\]
Unfortunately the above series is ``divergent''. We interpret it as the following convergent infinite sum of residues:
\[
\frac{1}{2^{12}} \sum_{n=1}^{\infty} {\rm Res} \left( \frac{(1)_{-s}^9 \Gamma(s)}{\left( \frac12 \right)_{-s}^7 \left( \frac14 \right)_{-s} \left( \frac34 \right)_{-s}} \frac{43680s^4+20632s^3+4340s^2+466s+21}{s^9} \frac{\cos \pi s}{2^{12s}} \right)_{s=n},
\]
which is $\overset?=$ to $95232 \zeta(5)$. Curiously the summand $-160 \pi^5 i$ does not appear when we include the factor $\cos \pi s$,

\subsection*{Conclusion}

Perhaps a kind of ``modular equations" (still undiscovered), combined with functional relations like those conjectured in Conj. \ref{conj2} and Conj. \ref{conj3}, could explain the algebraic values conjectured in Conj. \ref{conj1}.

\end{document}